\documentclass{article}
\usepackage{amssymb, amsmath}

\begin{document}

\begin{center}
{\Large A Clifford algebra associated to generalized Fibonacci quaternions}%
\begin{equation*}
\end{equation*}

Cristina FLAUT%
\begin{equation*}
\end{equation*}
\end{center}

\textbf{Abstract.} {\small In this paper we find a Clifford algebra
associated to generalized Fibonacci quaternions. In this way, we provide a
nice algorithm to obtain a division quaternion algebra starting from a
quaternion non-division algebra and vice-versa.} 
\begin{equation*}
\end{equation*}



\begin{equation*}
\end{equation*}

\textbf{1. Introduction}%
\begin{equation*}
\end{equation*}

In 1878, W. K. Clifford discovered Clifford algebras. These algebras
generalize the real numbers, complex numbers and quaternions( see [Le; 06 ]).

The theory of Clifford algebras is intimately connected with the theory of
quadratic forms. In the following, we will consider $K$ to be a field of
characteristic not two. Let $(V,q)$ be a $K-$vector space equipped with a
nondegenerate quadratic form over the field $K.$ A \textit{Clifford algebra }%
for $(V,q)$ is a $K-$algebra $C$ with a linear map $i:V\rightarrow C$ \
satisfying the property 
\begin{equation*}
i\left( x\right) ^{2}=q\left( x\right) \cdot 1_{C},\forall x\in V,
\end{equation*}%
such that for any $K-$algebra $A$ and any $K$ linear map $\gamma
:V\rightarrow A$ with $\gamma ^{2}\left( x\right) =q\left( x\right) \cdot
1_{A},\forall x\in V,$ there exists a unique \ $K$-algebra morphism $\gamma
^{\prime }:C\rightarrow A$ with $\gamma =\gamma ^{\prime }\circ i.$

Such an algebra can be constructed \ using the tensor algebra associated to
a vector space $V$. Let $T(V)=K\oplus V\oplus (V\otimes V)\oplus ...$ be the
tensor algebra associated to the vector space $V$ \ and let $\mathcal{J}$ be
the two-sided ideal of $T(V)$ generated by all elements \ of the form $%
x\otimes x-q\left( x\right) \cdot 1,$ for all $x\in V.$ The associated
Clifford \ algebra is the factor algebra $C(V,q)=T\left( V\right) /\mathcal{J%
}$ . \medskip ([Kn; 88], [La; 04])

\textbf{Theorem 1.1. (Poincar\'{e}-Birkhoff-Witt). (}[Kn; 88], p. 44\textbf{)%
}\textit{If} $\{e_{1},e_{2},...,e_{n}\}$ \textit{is a basis of} $V$ \textit{%
, then the set} $\{1,e_{j_{1}}e_{j_{2}}...e_{j_{s}},1\leq s\leq n,1\leq
j_{1}<j_{2}<...<j_{s}\leq n\}$ \textit{is a basis in} $C(V,q).\medskip $

The most important Clifford algebras are those defined over real and complex
vector spaces equipped with nondegenerate quadratic forms. Every
nondegenerate quadratic form over a real vector space is equivalent with the
following standard diagonal form:

\begin{equation*}
q(x)=x_{1}^{2}+...+x_{r}^{2}-x_{r+1}^{2}-...-x_{s}^{2},
\end{equation*}%
where $n=r+s$ is the dimension of the vector space. The pair of integers $%
(r,s)$ is called \textit{the signature} of the quadratic form. The real
vector space with this quadratic form is usually denoted $\mathbb{R}_{r,s}$
and the Clifford algebra on $\mathbb{R}_{r,s}$ is denoted $Cl_{r,s}\left( 
\mathbb{R}\right) $. For other details about Clifford algebras, the reader
is referred to [Ki, Ou; 99], [Ko; 10] and [Sm; 91].\medskip

\textbf{Example 1.2.}

i) For $p=q=0$ we have $Cl_{0,0}\left( K\right) \simeq K;$

ii) For $p=0,q=1,$ it results that $Cl_{0,1}\left( K\right) $ \ is a
two-dimensional algebra generated by a single vector $e_{1}$ such that $%
e_{1}^{2}=-1$ and therefore $Cl_{0,1}\left( K\right) \simeq K\left(
e_{1}\right) $. For $K=\mathbb{R}$ it follows that $Cl_{0,1}\left( \mathbb{R}%
\right) \simeq \mathbb{C}.$

iii) For $p=0,q=2,$ the algebra $Cl_{0,2}\left( K\right) $ is a
four-dimensional algebra spanned by the set $\{1,e_{1},e_{2},e_{1}e_{2}\}.$
Since $e_{1}^{2}=e_{2}^{2}=(e_{1}e_{2})^{2}=-1$ and $e_{1}e_{2}=-e_{2}e_{1},$
we obtain that this algebra is isomorphic to the division quaternions
algebra $\mathbb{H}$.

iv) For $p=1,q=1$ or $p=2,q=0,$ we obtain the algebra $Cl_{1,1}\left(
K\right) \simeq Cl_{2,0}\left( K\right) $ which is isomorphic with a
split(i.e. nondivision) quaternion algebra, called \textit{paraquaternion
algebra }or \textit{antiquaternion algebra}. (See [Iv, Za; 05])

\begin{equation*}
\end{equation*}

\textbf{2. Preliminaries}%
\begin{equation*}
\end{equation*}

Let $\mathbb{H}\left( \beta _{1},\beta _{2}\right) $ be the generalized
real\ quaternion algebra, the algebra of the elements of the form $%
a=a_{1}\cdot 1+a_{2}e_{2}+a_{3}e_{3}+a_{4}e_{4},$ where $a_{i}\in \mathbb{R}%
,i\in \{1,2,3,4\}$, and the elements of the basis $\{1,e_{2},e_{3},e_{4}\}$
satisfy the following multiplication table: \vspace{3mm}

\begin{center}
\begin{tabular}{ccccc}
$\cdot\:\,\,$\vline & $1$ & $e_2$ & $e_3$ & $e_4$ \\ \hline
$1$\:\,\vline & $1$ & $e_2$ & $e_3$ & $e_4$ \\ 
$e_2$\vline & $e_2$ & $-\beta_1$ & $e_4$ & $-\beta_1e_3$ \\ 
$e_3$\vline & $e_3$ & $-e_4$ & $-\beta_2$ & $\beta_2e_2$ \\ 
$e_4$\vline & $e_4$ & $\beta_1e_3$ & $-\beta_2e_2$ & $-\beta_1\beta_2$%
\end{tabular}
\end{center}

\vspace{3mm}

We denote by  $\boldsymbol{n}\left( a\right) $ the  norm of a real
quaternion $a.$ The norm of a generalized quaternion has the following
expression $\boldsymbol{n}\left( a\right) =a_{1}^{2}+\beta
_{1}a_{2}^{2}+\beta _{2}a_{3}^{2}+\beta _{1}\beta _{2}a_{4}^{2}.$ For $\beta
_{1}=\beta _{2}=1,$ we obtain the real division algebra $\mathbb{H},$ with
the basis $\{1,i,j,k\},$ where $i^{2}=j^{2}=k^{2}=-1$ and $%
ij=-ji,ik=-ki,jk=-kj$.\medskip 

\textbf{Proposition 2.1.} ([La; 04], Proposition 1.1) \textit{The quaternion
algebra} $\mathbb{H}\left( \beta _{1},\beta _{2}\right) $ \textit{is
isomorphic with quaternion algebra} $\mathbb{H}\left( x^{2}\beta
_{1},y^{2}\beta _{2}\right) ,$ \textit{where} $x,y\in K^{\ast }.$\textit{\ }$%
\Box \medskip $

The Fibonacci \ numbers was introduced \ by \textit{Leonardo of Pisa
(1170-1240) }in his book \textit{Liber abbaci}, book published in 1202 AD
(see [Kos; 01], p. 1, 3). This name is attached to the following sequence of
numbers%
\begin{equation*}
0,1,1,2,3,5,8,13,21,....,
\end{equation*}%
with the $n$th term given by the formula:%
\begin{equation*}
f_{n}=f_{n-1}+f_{n-2,}\ n\geq 2,\ 
\end{equation*}%
where $f_{0}=0,f_{1}=1.$

In [Ho; 61], the \ author generalized Fibonacci numbers and gave many
properties of them: 
\begin{equation*}
h_{n}=h_{n-1}+h_{n-2},\ \ n\geq 2,
\end{equation*}%
where $h_{0}=p,h_{1}=q,$ with $\ p,q$ being arbitrary integers. In the same
paper [Ho; 61], relation (7), the following relation between Fibonacci
numbers and generalized Fibonacci numbers was obtained: 
\begin{equation}
h_{n+1}=pf_{n}+qf_{n+1}.  \tag{2.1}
\end{equation}%
\qquad\ \ For the generalized real\ quaternion algebra, the Fibonacci
quaternions and generalized Fibonacci quaternions are defined in the same
way:%
\begin{equation*}
F_{n}=f_{n}\cdot 1+f_{n+1}e_{2}+f_{n+2}e_{3}+f_{n+3}e_{4},
\end{equation*}%
for the \ $n$th Fibonacci quaternions, and 
\begin{equation}
H_{n}=h_{n}\cdot 1+h_{n+1}e_{2}+h_{n+2}e_{3}+h_{n+3}e_{4}=pF_{n}+qF_{n+1}, 
\tag{2.2}
\end{equation}%
for the \ $n$th generalized \ Fibonacci quaternions.$\ $

In the following, we will denote the $n$th generalized Fibonacci number \
and a $n$th generalized Fibonacci quaternion element with $h_{n}^{p,q},$
respectively $H_{n}^{p,q}.$ In this way, we emphasis the starting integers $%
p $ and $q.\medskip $

It is known that the expression for the $n$th term of a Fibonacci element is 
\begin{equation}
f_{n}=\frac{1}{\sqrt{5}}[\alpha ^{n}-\beta ^{n}]=\frac{\alpha ^{n}}{\sqrt{5}}%
[1-\frac{\beta ^{n}}{\alpha ^{n}}],  \tag{2.3}
\end{equation}%
where $\alpha =\frac{1+\sqrt{5}}{2}$ and $\beta =\frac{1-\sqrt{5}}{2}.$

From the above, we obtain the following limit:

$\underset{n\rightarrow \infty }{\lim }\boldsymbol{n}\left( F_{n}\right) =%
\underset{n\rightarrow \infty }{\lim }(f_{n}^{2}+\beta _{1}f_{n+1}^{2}+\beta
_{2}f_{n+2}^{2}+\beta _{1}\beta _{2}f_{n+3}^{2})=$\newline
$=\underset{n\rightarrow \infty }{\lim }(\frac{\alpha ^{2n}}{5}$+$\beta _{1}%
\frac{\alpha ^{2n+2}}{5}$+$\beta _{2}\frac{\alpha ^{2n+4}}{5}$+$\beta
_{1}\beta _{2}\frac{\alpha ^{2n+6}}{5})=$\newline
$=sgnE(\beta _{1},\beta _{2})\cdot \infty ,$ where $E(\beta _{1},\beta _{2})=%
\frac{1}{5}[1+\beta _{1}+2\beta _{2}+5\beta _{1}\beta _{2}+\alpha \left(
\beta _{1}+3\beta _{2}+8\beta _{1}\beta _{2}\right) ],$ since $\alpha
^{2}=\alpha +1.$(see [Fl, Sh; 13])

If \ $E(\beta _{1},\beta _{2})>0,$ there exist a number $n_{1}\in \mathbb{N}$
such that for all\newline
$n\geq n_{1}$ we have $\boldsymbol{n}\left( F_{n}\right) >0.~$In the same
way, if $E(\beta _{1},\beta _{2})<0,$ there exist a number $n_{2}\in \mathbb{%
N}$ such that for all $n\geq n_{2}$ we have $\boldsymbol{n}\left(
F_{n}\right) <0.$ Therefore for \ all $\beta _{1},\beta _{2}\in \mathbb{R}$
with $E(\beta _{1},\beta _{2})\neq 0,$ in the algebra $\mathbb{H}\left(
\beta _{1},\beta _{2}\right) $ there is a natural number $n_{0}=\max
\{n_{1},n_{2}\}$ such that $\boldsymbol{n}\left( F_{n}\right) \neq 0,$ hence 
$F_{n}$ is an invertible element for all $n\geq n_{0}.$ Using the same
arguments, we can compute the following limit:\newline
$\underset{n\rightarrow \infty }{\lim }\left( \boldsymbol{n}\left(
H_{n}^{p,q}\right) \right) =\underset{n\rightarrow \infty }{\lim }\left(
h_{n}^{2}+\beta _{1}h_{n+1}^{2}+\beta _{2}h_{n+2}^{2}+\beta _{1}\beta
_{2}h_{n+3}^{2}\right) =sgnE^{\prime }(\beta _{1},\beta _{2})\cdot \infty ,$
where $E^{\prime }(\beta _{1},\beta _{2})=\frac{1}{5}\left( p+\alpha
q\right) ^{2}E(\beta _{1},\beta _{2}),$ if $E^{\prime }(\beta _{1},\beta
_{2})\neq 0.$(see [Fl, Sh; 13])

Therefore, for \ all $\beta _{1},\beta _{2}\in \mathbb{R}$ with $E^{\prime
}(\beta _{1},\beta _{2})\neq 0,$ in the algebra $\mathbb{H}\left( \beta
_{1},\beta _{2}\right) $ there exist a natural number $n_{0}^{\prime }$ such
that $\boldsymbol{n}\left( H_{n}^{p,q}\right) \neq 0,$ hence $H_{n}^{p,q}$
is an invertible element for all $n\geq n_{0}^{\prime }.\medskip $

\textbf{Theorem 2.2.} ([Fl, Sh; 13], Theorem 2.6 ) \textit{For all} $\beta
_{1},\beta _{2}\in \mathbb{R}$ \textit{with} $E^{\prime }(\beta _{1},\beta
_{2})\neq 0,$ \textit{there exist a natural number} $n^{\prime }$ \textit{%
such that for all} $n\geq n^{\prime }$ \textit{\ Fibonacci elements} $\
F_{n} $ \textit{and generalized Fibonacci elements} $H_{n}^{p,q}$ \textit{%
are invertible elements in the algebra} $\mathbb{H}\left( \beta _{1},\beta
_{2}\right) .\Box \medskip $

\textbf{Theorem} \textbf{2.3. }([Fl, Sh; 13], Theorem 2.1 )\textbf{\ } 
\textit{The set} $\mathcal{H}_{n}=\{H_{n}^{p,q}~/~p,q\in \mathbb{Z},n\geq
m,m\in \mathbb{N}\}\cup \{0\}$ \textit{is a} $\mathbb{Z}-$\textit{module}%
.\medskip $\Box \medskip $%
\begin{equation*}
\end{equation*}%
\textbf{3. Main results}%
\begin{equation*}
\end{equation*}

\textbf{Remark 3.1.} We remark that the $\mathbb{Z}-$module from \ Theorem
2.3 is a free $\mathbb{Z}-$module of rank $2$. Indeed, $\varphi :\mathbb{Z}%
\times \mathbb{Z\rightarrow }\mathcal{H}_{n},~\varphi \left( \left(
p,q\right) \right) =$ $H_{n}^{p,q}$ \ is a $\mathbb{Z}-$module isomorphism
and $\{$ $\varphi \left( 1,0\right) =F_{n},\varphi \left( 0,1\right)
=F_{n+1}\}$ is a basis in $\mathcal{H}_{n}$.\medskip

\textbf{Remark 3.2.} By extension of scalars, we obtain that $\mathbb{%
R\otimes }_{\mathbb{Z}}\mathcal{H}_{n}$ is a $\mathbb{R-}$vector space of
dimension two. A basis is $\{\overline{e}_{1}=1\mathbb{\otimes }F_{n},%
\overline{e}_{2}=1\mathbb{\otimes }F_{n+1}\}.$ We have that $\mathbb{%
R\otimes }_{\mathbb{Z}}\mathcal{H}_{n}$ is isomorphic with the $\mathbb{R-}$%
vector space $\mathcal{H}_{n}^{\mathbb{R}}=\{H_{n}^{p,q}~/~p,q\in \mathbb{R}%
\}\cup \{0\}.$ A basis in $\mathcal{H}_{n}^{\mathbb{R}}$ is $%
\{F_{n},F_{n+1}\}.$\medskip

Let $T\left( \mathcal{H}_{n}^{\mathbb{R}}\right) $ be the tensor algebra
associated to the $\mathbb{R-}$vector space $\mathcal{H}_{n}^{\mathbb{R}}$
and let $C\left( \mathcal{H}_{n}^{\mathbb{R}}\right) $ be the Clifford
algebra associated to tensor algebra $T\left( \mathcal{H}_{n}^{\mathbb{R}%
}\right) .$ From Theorem 1.1,  it results that this algebra has dimension
four.

\begin{equation*}
\end{equation*}

\textbf{Case} \textbf{1: }$\mathbb{H}\left( \beta _{1},\beta _{2}\right) $ 
\textbf{is a division algebra}%
\begin{equation*}
\end{equation*}

\textbf{Remark 3.3.} Since in this case\ $E(\beta _{1},\beta _{2})>0,$ for
all $n\geq n^{\prime }$ (as in Theorem 2.2.)$,$ then $\mathcal{H}_{n}^{%
\mathbb{R}}$ is an Euclidean vector space. Indeed, let $z,w\in \mathcal{H}%
_{n}^{\mathbb{R}%
},z=x_{1}F_{n}+x_{2}F_{n+1},w=y_{1}F_{n}+y_{2}F_{n+1},x_{1},x_{2},y_{1},y_{2}\in 
\mathbb{R}.$ The inner product is defined as in the following:%
\begin{equation*}
<z,w>=x_{1}y_{1}\mathbf{n}\left( F_{n}\right) +x_{2}y_{2}\mathbf{n}\left(
F_{n+1}\right) .
\end{equation*}%
We remark that all properties of inner product are fulfilled. Indeed, since
for all $n\geq n^{\prime }$ we have $\mathbf{n}\left( F_{n}\right) >0$ and $%
\mathbf{n}\left( F_{n+1}\right) >0,$ it results that $<z,z>=x_{1}^{2}\mathbf{%
n}\left( F_{n}\right) +x_{2}^{2}\mathbf{n}\left( F_{n+1}\right) =0$ if and
only if $x_{1}=x_{2}=0,$ therefore $z=0.$\medskip

On $\mathcal{H}_{n}^{\mathbb{R}}$ with the basis $\{F_{n},F_{n+1}\},~$we
define the following quadratic form $q_{\mathcal{H}_{n}^{\mathbb{R}}}:%
\mathcal{H}_{n}^{\mathbb{R}}\rightarrow \mathbb{R},$%
\begin{equation*}
q_{\mathcal{H}_{n}^{\mathbb{R}}}\left( x_{1}F_{n}+x_{2}F_{n+1}\right) =%
\mathbf{n}(F_{n})x_{1}^{2}+\mathbf{n}(F_{n+1})x_{2}^{2}.
\end{equation*}%
Let $Q_{\mathcal{H}_{n}^{\mathbb{R}}}$ be the bilinear form associated to
the quadratic form $q_{\mathcal{H}_{n}^{\mathbb{R}}},$%
\begin{eqnarray*}
Q_{\mathcal{H}_{n}^{\mathbb{R}}}\left( x,y\right) &=&\frac{1}{2}(q_{\mathcal{%
H}_{n}^{\mathbb{R}}}(x+y)-q_{\mathcal{H}_{n}^{\mathbb{R}}}(x)-q_{\mathcal{H}%
_{n}^{\mathbb{R}}}(y))= \\
&=&\mathbf{n}(F_{n})x_{1}y_{1}+\mathbf{n}(F_{n+1})x_{2}y_{2}.
\end{eqnarray*}%
The matrix associated to the quadratic form $q_{\mathcal{H}_{n}^{\mathbb{R}%
}} $ is 
\begin{equation*}
A=\left( 
\begin{array}{cc}
\mathbf{n}(F_{n}) & 0 \\ 
0 & \mathbf{n}(F_{n+1})%
\end{array}%
\right) .
\end{equation*}%
\newline

We remark that $\det A=\mathbf{n}(F_{n})\mathbf{n}(F_{n+1})>0,$ for all $%
n\geq n^{\prime }$. Since $E(\beta _{1},\beta _{2})>0,$ therefore $\mathbf{n}%
(F_{n})>0,$ for $n>n^{\prime }.$ We obtain that the quadratic form $q_{%
\mathcal{H}_{n}^{\mathbb{R}}}$ is positive definite and the Clifford algebra 
$C\left( \mathcal{H}_{n}^{\mathbb{R}}\right) $ associated to the tensor
algebra $T\left( \mathcal{H}_{n}^{\mathbb{R}}\right) $ is isomorphic with $%
Cl_{2,0}\left( K\right) $ which is isomorphic to a split quaternion algebra.

From the above results \ and using Proposition 2.1, we obtain the following
theorem:\medskip 

\textbf{Theorem 3.4.} \textit{If } $\mathbb{H}\left( \beta _{1},\beta
_{2}\right) $ \textit{is a division algebra, there is a natural number} $%
n^{\prime }$ \textit{such that} \textit{for all} $n\geq n^{\prime },$ 
\textit{the Clifford algebra associated to the real vector space} $\mathcal{H%
}_{n}^{\mathbb{R}}$ \textit{is isomorphic with the split quaternion algebra} 
$\mathbb{H}\left( -1,-1\right) .$\textit{\ }$\Box \medskip $

\begin{equation*}
\end{equation*}

\bigskip \textbf{Case} \textbf{2: } $\mathbb{H}\left( \beta _{1},\beta
_{2}\right) $ \textbf{is not a division algebra}

\begin{equation*}
\end{equation*}

\textbf{Remark 3.5.} i) If \ $E(\beta _{1},\beta _{2})>0,$ then $\mathcal{H}%
_{n}^{\mathbb{R}}$ is an Euclidean vector space, for all $n\geq n^{\prime }$%
, as in Theorem 2.2. Indeed, let $z,w\in $ $\mathcal{H}_{n}^{\mathbb{R}%
},z=x_{1}F_{n}+x_{2}F_{n+1},w=y_{1}F_{n}+y_{2}F_{n+1},x_{1},x_{2},y_{1},y_{2}\in 
\mathbb{R}.$ The inner product is defined as in the following:%
\begin{equation*}
<z,w>=x_{1}y_{1}\mathbf{n}\left( F_{n}\right) +x_{2}y_{2}\mathbf{n}\left(
F_{n+1}\right) .
\end{equation*}

ii) If \ $E(\beta _{1},\beta _{2})<0,$ then $\mathcal{H}_{n}^{\mathbb{R}}$
is also an Euclidean vector space, for all $n\geq n^{\prime }$, as in
Theorem 2.2. \ Indeed, let $z,w\in $ $\mathcal{H}_{n}^{\mathbb{R}%
},z=x_{1}F_{n}+x_{2}F_{n+1},w=y_{1}F_{n}+y_{2}F_{n+1},x_{1},x_{2},y_{1},y_{2}\in 
\mathbb{R}.$ The inner product is defined as in the following:%
\begin{equation*}
<z,w>=-x_{1}y_{1}\mathbf{n}\left( F_{n}\right) -x_{2}y_{2}\mathbf{n}\left(
F_{n+1}\right) .
\end{equation*}%
We have $<z,z>=-x_{1}^{2}\mathbf{n}\left( F_{n}\right) -x_{2}^{2}\mathbf{n}%
\left( F_{n+1}\right) ,$ and since for all $n\geq n^{\prime }$ we have $%
\mathbf{n}\left( F_{n}\right) <0$ and $\mathbf{n}\left( F_{n+1}\right) <0,$
it results that $<z,z>=-x_{1}^{2}\mathbf{n}\left( F_{n}\right) -x_{2}^{2}%
\mathbf{n}\left( F_{n+1}\right) =0$ if and only if $x_{1}=x_{2}=0,$
therefore $z=0.\medskip $

On $\mathcal{H}_{n}^{\mathbb{R}}$ with the basis $\{F_{n},F_{n+1}\},~$we
define the following quadratic form $q_{\mathcal{H}_{n}^{\mathbb{R}}}:%
\mathcal{H}_{n}^{\mathbb{R}}\rightarrow \mathbb{R},$%
\begin{equation*}
q_{\mathcal{H}_{n}^{\mathbb{R}}}\left( x_{1}F_{n}\text{+}x_{2}F_{n+1}\right)
=q_{\mathcal{H}_{n}^{\mathbb{R}}}\left( x_{1}F_{n}+x_{2}F_{n+1}\right) =%
\mathbf{n}\left( F_{n}\right) x_{1}^{2}+\mathbf{n}\left( F_{n+1}\right)
x_{2}^{2}.
\end{equation*}%
Let $Q_{\mathcal{H}_{n}^{\mathbb{R}}}$ be the bilinear form associated to
the quadratic form $q_{\mathcal{H}_{n}^{\mathbb{R}}},$%
\begin{eqnarray*}
Q_{\mathcal{H}_{n}^{\mathbb{R}}}\left( x,y\right) &=&\frac{1}{2}(q_{\mathcal{%
H}_{n}^{\mathbb{R}}}(x+y)-q_{\mathcal{H}_{n}^{\mathbb{R}}}(x)-q_{\mathcal{H}%
_{n}^{\mathbb{R}}}(y))= \\
&=&\mathbf{n}\left( F_{n}\right) x_{1}y_{1}+\mathbf{n}\left( F_{n+1}\right)
x_{2}y_{2}.
\end{eqnarray*}%
The matrix associated to quadratic form $q_{\mathcal{H}_{n}^{\mathbb{R}}}$
is 
\begin{equation*}
A=\left( 
\begin{array}{cc}
\mathbf{n}(F_{n}) & 0 \\ 
0 & \mathbf{n}(F_{n+1})%
\end{array}%
\right) .
\end{equation*}%
\newline

We remark that $\det A=\mathbf{n}(F_{n})\mathbf{n}(F_{n+1})>0,$ for all $%
n\geq n^{\prime }$.

If $E(\beta _{1},\beta _{2})>0,$ therefore $\mathbf{n}(F_{n})>0,$ for $%
n>n^{\prime }.$ We obtain that the quadratic form $q_{\mathcal{H}_{n}^{%
\mathbb{R}}}$ is positive definite and the Clifford algebra $C\left( 
\mathcal{H}_{n}^{\mathbb{R}}\right) $ associated to the tensor algebra $%
T\left( \mathcal{H}_{n}^{\mathbb{R}}\right) $ is isomorphic with $%
Cl_{2,0}\left( K\right) $ which is isomorphic to a split quaternion algebra.

If $E(\beta _{1},\beta _{2})<0,$ therefore $\mathbf{n}(F_{n})<0,$ for $%
n>n^{\prime }.$ Then the quadratic form $q_{\mathcal{H}_{n}^{\mathbb{R}}}$
is negative definite and the Clifford algebra $C\left( \mathcal{H}_{n}^{%
\mathbb{R}}\right) $ associated to the tensor algebra $T\left( \mathcal{H}%
_{n}^{\mathbb{R}}\right) $ is isomorphic with $Cl_{0,2}\left( K\right) $
which is isomorphic to the quaternion division algebra $\mathbb{H}$.

From the above results \ and using Proposition 2.1, we obtain the following
theorem:\medskip 

\textbf{Theorem 3.6.} \textit{If } $\mathbb{H}\left( \beta _{1},\beta
_{2}\right) $ \textit{is not a division algebra, there is a natural number} $%
n^{\prime }$ \textit{such that} \textit{for all} $n\geq n^{\prime },$ 
\textit{if }\ $E(\beta _{1},\beta _{2})>0,$ \textit{then the Clifford
algebra associated to the real vector space} $\mathcal{H}_{n}^{\mathbb{R}}$ 
\textit{is isomorphic with the split quaternion algebra} $\mathbb{H}\left(
-1,-1\right) .$ If $E(\beta _{1},\beta _{2})<0,$ \textit{then the Clifford
algebra associated to the real vector space} $\mathcal{H}_{n}^{\mathbb{R}}$ 
\textit{is isomorphic with the division quaternion algebra} $\mathbb{H}%
\left( 1,1\right) .$\textit{\ }$\Box \medskip $

\textbf{Example 3.7.} 1) For \ $\beta _{1}=1,\beta _{2}=-1,$ we obtain the
split quaternion algebra $\mathbb{H}\left( 1,-1\right) $. In this case, we
have $E(\beta _{1},\beta _{2})=\frac{1}{5}[-5-10\alpha ]<0$ and, for $%
n^{\prime }=0,$ we obtain $\mathbf{n}\left( F_{n}\right) =$ $%
f_{n}^{2}+f_{n+1}^{2}-f_{n+2}^{2}-f_{n+3}^{2}<0,\mathbf{n}\left(
F_{n+1}\right) =$ $f_{n+1}^{2}+f_{n+2}^{2}-f_{n+3}^{2}-f_{n+4}^{2}<0,$ for
all $n\geq 0.$ The quadratic form $q_{\mathcal{H}_{n}^{\mathbb{R}}}$ is
negative definite, therefore the Clifford algebra $C\left( \mathcal{H}_{n}^{%
\mathbb{R}}\right) $ associated to the tensor algebra $T\left( \mathcal{H}%
_{n}^{\mathbb{R}}\right) $ is isomorphic to $Cl_{0,2}\left( K\right) $ which
is isomorphic to the quaternion division algebra $\mathbb{H}\left(
1,1\right) .$

2) For $\beta _{1}=-2,\beta _{2}=-3,$ we obtain the split quaternion algebra 
$\mathbb{H}\left( -2,-3\right) .$ In this case, we have $E(\beta _{1},\beta
_{2})=\frac{1}{5}\left[ 23+43\alpha \right] >0.$ For $n^{\prime }=0,$ we
obtain $\mathbf{n}\left( F_{n}\right) =$ $%
f_{n}^{2}-f_{n+1}^{2}-f_{n+2}^{2}+f_{n+3}^{2}>0,\mathbf{n}\left(
F_{n+1}\right) =$ $f_{n+1}^{2}-f_{n+2}^{2}-f_{n+3}^{2}+f_{n+4}^{2}>0,$ for
all $n\geq 0.$ The quadratic form $q_{\mathcal{H}_{n}^{\mathbb{R}}}$ is
positive definite, therefore the Clifford algebra $C\left( \mathcal{H}_{n}^{%
\mathbb{R}}\right) $ associated to the tensor algebra $T\left( \mathcal{H}%
_{n}^{\mathbb{R}}\right) $ is isomorphic to $Cl_{2,0}\left( K\right) $ which
is isomorphic to the split quaternion algebra $\mathbb{H}\left( -1,-1\right)
.$

3) For $\beta _{1}=2,\beta _{2}=-3,$ we obtain the split quaternion algebra $%
\mathbb{H}\left( 2,-3\right) .$ In this case, we have $E(\beta _{1},\beta
_{2})=\frac{1}{5}\left[ -33-44\alpha \right] <0.$ For $n^{\prime }=0,$ we
obtain $\mathbf{n}\left( F_{n}\right) =$ $%
f_{n}^{2}+2f_{n+1}^{2}-3f_{n+2}^{2}-6f_{n+3}^{2}<0,\mathbf{n}\left(
F_{n+1}\right) =$ $f_{n+1}^{2}+2f_{n+2}^{2}-3f_{n+3}^{2}-6f_{n+4}^{2}>0,$
for all $n\geq 0.$ The quadratic form $q_{\mathcal{H}_{n}^{\mathbb{R}}}$ is
negative definite, therefore the Clifford algebra $C\left( \mathcal{H}_{n}^{%
\mathbb{R}}\right) $ associated to the tensor algebra $T\left( \mathcal{H}%
_{n}^{\mathbb{R}}\right) $ is isomorphic to $Cl_{0,2}\left( K\right) $ which
is isomorphic to the division quaternion algebra $\mathbb{H}\left(
1,-1\right) .$

3) For $\beta _{1}=\beta _{2}=-\frac{1}{2},$ we obtain the split quaternion
algebra $\mathbb{H}\left( -\frac{1}{2},-\frac{1}{2}\right) .$ Therefore $%
E(\beta _{1},\beta _{2})=\frac{3}{20}>0$ and for $n^{\prime }=1$ we obtain $%
\mathbf{n}\left( F_{n}\right) >0$ and $\mathbf{n}\left( F_{n+1}\right) >0.$%
The quadratic form $q_{\mathcal{H}_{n}^{\mathbb{R}}}$ is positive definite,
therefore the Clifford algebra $C\left( \mathcal{H}_{n}^{\mathbb{R}}\right) $
associated to the tensor algebra $T\left( \mathcal{H}_{n}^{\mathbb{R}%
}\right) $ is isomorphic with $Cl_{2,0}\left( K\right) $ which is isomorphic
to the split quaternion algebra $\mathbb{H}\left( -1,-1\right) .\medskip $

\textbf{The algorithm.}

1) Let \ $\mathbb{H}\left( \beta _{1},\beta _{2}\right) $ be a quaternion
algebra, $\alpha =\frac{1+\sqrt{5}}{2}$ and $E(\beta _{1},\beta _{2})=\frac{1%
}{5}[1+\beta _{1}+2\beta _{2}+5\beta _{1}\beta _{2}+\alpha \left( \beta
_{1}+3\beta _{2}+8\beta _{1}\beta _{2}\right) ],$

2) Let $V$ be the $\mathbb{R-}$vector space $\mathcal{H}_{n}^{\mathbb{R}%
}=\{H_{n}^{p,q}~/~p,q\in \mathbb{R}\}\cup \{0\}.$

3) If $E(\beta _{1},\beta _{2})>0,$ then the Clifford algebra $C\left( 
\mathcal{H}_{n}^{\mathbb{R}}\right) $ associated to the tensor algebra $%
T\left( \mathcal{H}_{n}^{\mathbb{R}}\right) $ is isomorphic with $%
Cl_{2,0}\left( K\right) $ which is isomorphic to the split quaternion
algebra $\mathbb{H}\left( -1,-1\right) .$

4) If $E(\beta _{1},\beta _{2})<0,$ then the Clifford algebra $C\left( 
\mathcal{H}_{n}^{\mathbb{R}}\right) $ associated to the tensor algebra $%
T\left( \mathcal{H}_{n}^{\mathbb{R}}\right) $ is isomorphic with $%
Cl_{0,2}\left( K\right) $ which is isomorphic to the division quaternion
algebra $\mathbb{H}\left( 1,1\right) .\medskip $

\textbf{Conclusions.} In this paper, we extend the $\mathbb{Z}-$module of
the generalized Fibonacci quaternions to a real vector space $\mathcal{H}%
_{n}^{\mathbb{R}}.$ We proved that the Clifford algebra $C\left( \mathcal{H}%
_{n}^{\mathbb{R}}\right) $ associated to the tensor algebra $T\left( 
\mathcal{H}_{n}^{\mathbb{R}}\right) $ is isomorphic to a split quaternion
algebra or to a division algebra if $E(\beta _{1},\beta _{2})=\frac{1}{5}%
[1+\beta _{1}+2\beta _{2}+5\beta _{1}\beta _{2}+\alpha \left( \beta
_{1}+3\beta _{2}+8\beta _{1}\beta _{2}\right) ]$ is positive or negative. We
also gave an algorithm which allows us to find a division quaternion algebra
starting from a split quaternion algebra and vice-versa.

\begin{equation*}
\end{equation*}

\textbf{References}\medskip

[Fl, Sh; 13] C. Flaut, V. Shpakivskyi, \textit{On Generalized Fibonacci
Quaternions and Fibonacci-Narayana Quaternions}, Adv. Appl. Clifford
Algebras, \textbf{23(3)}(2013), 673-688.

[Ho; 61] A. F. Horadam, \textit{A Generalized Fibonacci Sequence}, Amer.
Math. Monthly, \textbf{68}(1961), 455-459.

[Ki, Ou; 99] \ El Kinani, E. H., Ouarab, A., \textit{The Embedding of} $%
U_{q}(sl\left( 2\right) )$ \textit{and Sine Algebras in Generalized Clifford
Algebras}, Adv. Appl. Clifford Algebr., \textbf{9(1)}(1999), 103-108.

[Ko; 10] Ko\c{c}, C., \textit{C-lattices and decompositions of generalized
Clifford algebras,} Adv. Appl. Clifford Algebr., \textbf{20(2)}(2010),
313-320.

[Kos; 01] T.~Koshy, \textit{Fibonacci and Lucas Numbers with Applications},
A Wiley-Interscience publication, U.S.A, 2001.

[Kn; 88] M. A. Knus, \textit{Quadratic Forms, Clifford Algebras and Spinors}%
, IMECC-UNICAMP, 1988.

[La; 04] T.Y. Lam, \textit{Quadratic forms over fields}, AMS, Providence,
Rhode Island, 2004.

[Le; 06 ] Lewis, D. W., \textit{Quaternion Algebras and the Algebraic Legacy
of Hamilton's Quaternions}, Irish Math. Soc. Bulletin \textbf{57}(2006),
41--64.

[Sm; 91] Smith T. L., \ \textit{Decomposition of Generalized Clifford
Algebras}, Quart. J. Math. Oxford, \textbf{42}(1991), 105-112.

\begin{equation*}
\end{equation*}

Cristina FLAUT

{\small Faculty of Mathematics and Computer Science,}

{\small Ovidius University,}

{\small Bd. Mamaia 124, 900527, CONSTANTA,}

{\small ROMANIA}

{\small http://cristinaflaut.wikispaces.com/}

{\small http://www.univ-ovidius.ro/math/}

{\small e-mail:}

{\small cflaut@univ-ovidius.ro}

{\small cristina\_flaut@yahoo.com}%
\begin{equation*}
\end{equation*}

\end{document}